# EXISTENCE, UNIQUENESS AND $L_t^2(H_x^2) \cap L_t^\infty(H_x^1) \cap H_t^1(L_x^2)$ REGULARITY OF THE GRADIENT FLOW OF THE AMBROSIO–TORTORELLI FUNCTIONAL

## Tommaso Cortopassi

**Abstract.** We consider the gradient flow of the Ambrosio–Tortorelli functional, proving existence, uniqueness and $L_t^2(H_x^2) \cap L_t^\infty(H_x^1) \cap H_t^1(L_x^2)$ regularity of the solution in dimension 2. Such functional is an approximation in the sense of $\Gamma$-convergence of the Mumford–Shah functional often used in problems of image segmentation and fracture mechanics. The strategy of the proof essentially follows the one of []Feng and Prohl, *M2AN Math. Model. Numer. Anal.* **38** (2004) 291–320] but the crucial estimate is attained employing a different technique, and in the end it allows to prove better estimates than the ones obtained in [Feng and Prohl, *M2AN Math. Model. Numer. Anal.* **38** (2004) 291–320]. In particular we prove that if $U \subset \mathbb{R}^2$ is a bounded $C^2$ domain, the initial data $(u_0, z_0) \in [H^1(U)]^2$ with $0 \le z_0 \le 1$, then for every $T > 0$ there exists a unique gradient flow $(u(t), z(t))$ of the Ambrosio–Tortorelli functional such that

$$(u, z) \in [L^2(0, T; H^2(U)) \cap L^\infty(0, T; H^1(U)) \cap H^1(0, T; L^2(U))]^2,$$

while previously such regularity was known only for short times.

## 1. Introduction

The Mumford–Shah functional, introduced in [1], is defined as

$$E(u, \Gamma) = \frac{1}{2} \int_{U \setminus \Gamma} |\nabla u|^2 + (u - g)^2 \mathrm{d}x + \mathcal{H}^1(\Gamma), \tag{1.1}$$

where $u \in H^1(U \setminus \Gamma)$, $U \subset \mathbb{R}^2$ is open and bounded, $\Gamma \subset U$ is closed and $\mathcal{H}^1$ is the one dimensional Hausdorff measure. This functional has been extensively studied for its applications in image segmentation and fracture mechanics (a detailed survey can be found in [2]). The definition of the functional depends on the model it is used for: in (1.1) we gave the original definition of [1] which is suited for image segmentation models, while the interested reader may want to see the seminal paper [3] about the fracture mechanics case. The idea is rather simple: given a gray image $g$ (*i.e.* a scalar function), we want to find $(\tilde{u}, \tilde{\Gamma})$ such that

$$(\tilde{u}, \tilde{\Gamma}) = \operatorname*{arg\,min}_{\substack{u \in H^1(U \setminus \Gamma) \\ \Gamma \subset U \text{ closed}}} E(u, \Gamma). \tag{1.2}$$

In this case, the function $u$ will approximate in a "smooth way" the image $g$ while $\Gamma$ will be the contour set. From a theoretical point of view, in [4] the authors proved the existence of a solution for (1.2) by restricting to functions $u \in SBV(U)$ (see [5]) and $\Gamma = S_u$, *i.e.* the set of discontinuity jump points of $u$. From a numerical perspective, since the functional involves the measure of the singular jump set of a (unknown) function $u$, its direct numerical implementation is problematic. A standard approach is to minimize a more regular functional proposed by Ambrosio and Tortorelli in [6], [7] defined as

$$J_\varepsilon(u, z) = \frac{1}{2} \int_U [(\eta_\varepsilon + z^2)|\nabla u|^2 + (u - g)^2] \mathrm{d}x + \int_U \left[ \frac{(1-z)^2}{4\varepsilon} + \varepsilon |\nabla z|^2 \right] \mathrm{d}x \tag{1.3}$$

with $\eta_\varepsilon, \varepsilon > 0$, which approximates the Mumford–Shah functional in the sense of $\Gamma$-convergence as $\varepsilon \to 0$ in arbitrary dimension. The Ambrosio–Tortorelli functional provides a way which is usually employed to find (approximate) minima points for (1.1), *i.e.* to use a gradient flow approach for minimizing (1.3). As in [8] we will consider the gradient flow of (1.3), which is given by

$$\begin{cases} \partial_t u = \operatorname{div}((\eta_\varepsilon + z^2)\nabla u) - (u - g) & \text{in } (0, T) \times U \\ \partial_t z = 2\varepsilon \Delta z - z|\nabla u|^2 + \frac{1-z}{2\varepsilon} & \text{in } (0, T) \times U \\ \partial_n u = \partial_n z = 0 & \text{in } (0, T) \times \partial U \\ u(0, \cdot) = u_0 \text{ and } z(0, \cdot) = z_0 & \text{in } U \end{cases} \tag{1.4}$$

with $\eta_\varepsilon, \varepsilon > 0$ fixed. As already mentioned, our goal is to study the existence, uniqueness and regularity of solutions of (1.4), and the novelty lies in an approach (already suggested in a footnote of [8] and used in [9] for a different problem) which simplifies the proof of [8] while also gaining better regularity results. In particular in [8] they only manage to prove the $L_t^2(H_x^2)$ regularity for a short time $T_1 = T_1(\varepsilon)$ such that $T_1 \to 0$ as $\varepsilon \to 0$, with the crucial estimate being a local energy estimate inspired by a technique due to Struwe (see [10]) which only holds for sufficiently small times. We will prove $L_t^2(H_x^2)$ regularity for all positive times with the additional assumption that $g$ in (1.4) is not time dependent. The key estimate will be in Proposition 2.10, where we bound the $L_t^4(L_x^4)$ norms of the gradients $\nabla u, \nabla z$.

## 2. MAIN RESULT

Let $U \subset \mathbb{R}^2$ be a $C^2$ bounded domain and let $(u_0, z_0) \in [H^1(U)]^2$ with $0 \le z_0 \le 1$, $g \in L^2(U)$. As stated in the introduction, we want to prove existence, uniqueness and regularity of solutions of the gradient flow of (1.3), given by (1.4). However we will mostly work with the slightly modified system

$$\begin{cases} \partial_t u = \operatorname{div}((\eta_\varepsilon + \phi(z)^2)\nabla u) - (u - g) & \text{in } (0, T) \times U \\ \partial_t z = 2\varepsilon \Delta z - \phi'(z)\phi(z)|\nabla u|^2 + \frac{1-z}{2\varepsilon} & \text{in } (0, T) \times U \\ \partial_n u = \partial_n z = 0 & \text{in } (0, T) \times \partial U \\ u(0, \cdot) = u_0 \text{ and } z(0, \cdot) = z_0 & \text{in } U \end{cases} \tag{2.1}$$

where $\phi$ is a nondecreasing smooth cutoff function defined as

$$\phi(s) = \begin{cases} -1 & \text{if } s \leq -1 \\ s & \text{if } 0 < s < 1 \\ 2 & \text{if } s \geq 2. \end{cases}$$

The reason to introduce such a cutoff, as we will show later, is to bound the $L^\infty$ norm of $\eta_\varepsilon + \phi(z_N)^2$ when considering approximate solutions $z_N$ that we will define later. We will work with a modified Ambrosio–Tortorelli functional (of which (2.1) is the gradient flow of) denoted as $AT_\varepsilon$ and defined as

$$AT_\varepsilon(u(t), z(t)) = \frac{1}{2} \int_U [(\eta_\varepsilon + \phi(z(t))^2)|\nabla u(t)|^2 + (u(t) - g)^2] \mathrm{d}x \tag{2.2}$$
$$+ \int_U \left[ \frac{(1 - z(t))^2}{4\varepsilon} + \varepsilon |\nabla z(t)|^2 \right] \mathrm{d}x.$$

First of all, we define the notions of strong and weak solutions for (1.4) and (2.1).

**Definition 2.1. (Weak and strong solutions)**
We say that a couple $(u, z) \in [L^\infty(0, T; L^2(U)) \cap L^2(0, T; H^1(U))]^2$ is a weak solution to (1.4) if for almost every $t \in [0, T]$ it holds

$$\int_U u(t, x)\xi(x)\mathrm{d}x = \int_U u_0(x)\xi(x)\mathrm{d}x - \int_0^t \int_U (\eta_\varepsilon + z^2)\nabla u \nabla \xi \mathrm{d}x \mathrm{d}s$$
$$- \int_0^t \int_U (u - g)\xi \mathrm{d}x \mathrm{d}s \quad \forall \xi \in C_c^\infty(U).$$

and

$$\int_U z(t, x)\xi(x)\mathrm{d}x = \int_U z(0, x)\xi(x)\mathrm{d}x - \int_0^t \int_U 2\varepsilon \nabla z \nabla \xi \mathrm{d}x \mathrm{d}s$$
$$- \int_0^t \int_U z\xi|\nabla u|^2 \mathrm{d}x \mathrm{d}s + \int_0^t \int_U \frac{1 - z}{2\varepsilon} \xi \mathrm{d}x \mathrm{d}s \quad \forall \xi \in C_c^\infty(U).$$

If $(u, z) \in [L^2(0, T; H^2(U)) \cap L^\infty(0, T; H^1(U)) \cap H^1(0, T; L^2(U))]^2$ we say it is a strong solution. We define in the same way weak and strong solutions for (2.1).

Before starting to prove the necessary preliminary lemmas, let us state the main result of the paper.

**Theorem 2.2. *(Main result)***
*Let $U \subset \mathbb{R}^2$ be a bounded $C^2$ domain and let $(u_0, z_0) \in [H^1(U)]^2$ with $0 \leq z_0 \leq 1$ and $g \in L^2(U)$. Then there exists a unique strong solution $(u, z)$ of (1.4), with $\varepsilon, \eta_\varepsilon > 0$ fixed.*

**Remark 2.3.** It is worth to point out that all the estimates we will get are limited to $\eta_\varepsilon, \varepsilon > 0$ fixed. When $\varepsilon \to 0$ it may be that such bounds blow up.

We start with a proposition which shows how a strong solution of (2.1) is also a strong solution of (1.4) if $0 \leq z_0 \leq 1$.

**Proposition 2.4.** *Let $0 \leq z_0 \leq 1$. If a strong solution $(u, z)$ for (2.1) exists, it must be $0 \leq z(t, x) \leq 1$ for every $t \in (0, T)$ and for a.e. $x \in U$. So $(u, z)$ will also be a strong solution of (1.4).*

*Proof.* Consider

$$f_0(z) \coloneqq \max\{-z, 0\} \text{ and } f_1(z) \coloneqq \max\{z - 1, 0\}. \tag{2.3}$$

Differentiating $||f_0(z(t))||^2_{L^2(U)}$ in time we have by Reynolds transport theorem (we omit explicitly writing the dependence on $x$):

$$\frac{1}{2}\frac{\mathrm{d}}{\mathrm{d}t}||f_0(z(t))||^2_{L^2(U)} = \frac{\mathrm{d}}{\mathrm{d}t}\int_{\{z(t)<0\}}\frac{1}{2}z(t)^2\mathrm{d}x = \int_{\{z(t)<0\}}z(t)\partial_t z(t)\mathrm{d}x \tag{2.4}$$

$$+ \underbrace{\int_{\partial\{z(t)<0\}}\langle\omega(t),\hat{n}(t)\rangle\frac{1}{2}z(t)^2\mathrm{d}\mathcal{H}^1(x)}_{=0} = \int_U \chi_{\{z(t)<0\}}z(t)\partial_t z(t)\mathrm{d}x$$

with $\hat{n}$ the outward normal to $\partial\{z(t) < 0\}$, $\omega(t, x)$ the velocity of the boundary $\partial\{z(t) < 0\}$ and $\langle\cdot,\cdot\rangle$ the standard scalar product. Notice that the integral over $\partial\{z(t) < 0\}$ in (2.4) is 0 since $\langle\omega(t),\hat{n}(t)\rangle z(t) = 0$ on $\partial\{z(t) < 0\}$. Using that $z$ solves (2.1) by assumption:

$$\frac{1}{2}\frac{\mathrm{d}t}{\mathrm{d}t}||f_0(z(t))||^2_{L^2(U)} = \int_{\{z(t)<0\}}z(t)\partial_t z(t)\mathrm{d}x \tag{2.5}$$

$$= \int_{\{z(t)<0\}}2\varepsilon z(t)\Delta z(t) - z(t)\phi'(z(t))\phi(z(t))|\nabla u(t)|^2 + \frac{1-z(t)}{2\varepsilon}z(t)\mathrm{d}x$$

$$= -\int_{\{z(t)<0\}}2\varepsilon|\nabla z(t)|^2\mathrm{d}x + 2\varepsilon\underbrace{\int_{\{z(t)<0\}}\frac{\partial z}{\partial n}(t)z(t)\mathrm{d}\mathcal{H}^1(x)}_{=0}$$

$$- \int_{\{z(t)<0\}}z(t)\phi'(z(t))\phi(z(t))|\nabla u(t)|^2\mathrm{d}x + \int_{\{z(t)<0\}}\frac{1-z(t)}{2\varepsilon}z(t)\mathrm{d}x,$$

and each integral is nonpositive, except for the integral over $\partial\{z(t) < 0\}$ which is 0, since $z(t)\partial_n z(t) = 0$ on $\partial\{z(t) < 0\}$. Having that $f_0(z(0)) = 0$ because $z_0 \geq 0$, we have $f_0(z(t)) \equiv 0$ for every $t$, so $z(t) \geq 0$. In the same way we can prove $z(t) \leq 1$ by considering $f_1(z(t))$. $\square$

The strategy to prove Theorem 2.2 is to make use of Galerkin approximations, *i.e.* approximate solutions on finite dimensional spaces. We consider an orthogonal basis of $H^1(U)$ composed of eigenfunctions of $-\Delta$ on $U$ with homogeneous Neumann boundary conditions normalized with respect to the $L^2(U)$ norm, and denote it as $\{e_i\}_{i\in\mathbb{N}}$. So

$$\begin{cases}-\Delta e_i = \lambda_i e_i & \text{in } U \\ \partial_n e_i = 0 & \text{in } \partial U\end{cases} \text{ and } ||e_i||_{L^2(U)} = 1 \text{ for every } i. \tag{2.6}$$

We want to find Galerkin approximates $(u_N, z_N)$ such that they solve (in distributional sense) in $V_N = \mathrm{Span}(\{e_1, \ldots, e_N\})$ the system

$$\begin{cases}\partial_t u_N = \pi_N[\mathrm{div}((\eta_\varepsilon + \phi(z_N)^2)\nabla u_N)] - (u_N - g_N) & \text{in } (0, T) \times U \\ \partial_t z_N = 2\varepsilon\Delta z_N - \pi_N[\phi'(z_N)\phi(z_N)|\nabla u_N|^2] + \frac{1-z_N}{2\varepsilon} & \text{in } (0, T) \times U \\ \partial_n u_N = \partial_n z_N = 0 & \text{in } (0, T) \times \partial U \\ u_N(0, \cdot) = \pi_N[u_0] \text{ and } z_N(0, \cdot) = \pi_N[z_0] & \text{in } U\end{cases} \tag{2.7}$$

with $\pi_N$ the orthogonal projection on $V_N$. Notice that by the orthogonality properties of $\{e_i\}_{i=1}^N$, (2.7) is a $2N$ system of ODEs in $u^{(1)}, \ldots, u^{(N)}, z^{(1)}, \ldots, z^{(N)}$ with

$$u_N(t, x) = \sum_{i=1}^N u^{(i)}(t)e_i(x); \qquad z_N(t, x) = \sum_{i=1}^N z^{(i)}(t)e_i(x) \tag{2.8}$$

$$\pi_N[u_0] = \sum_{i=1}^N (u_0)^{(i)}e_i(x); \qquad \pi_N[z_0] = \sum_{i=1}^N (z_0)^{(i)}e_i(x). \tag{2.9}$$

which by the standard Cauchy–Lipschitz theorem admits a unique local solution. Indeed if we test (2.7) with $e_i$ we get:

$$u^{(i)}(t) = (u_0)^{(i)} - \int_0^t \int_U (\eta_\varepsilon + \phi(z_N(s))^2)\nabla u_N(s)\nabla e_i + (u_N(s) - g_N)e_i \mathrm{d}x\mathrm{d}s \tag{2.10}$$

$$z^{(i)}(t) = (z_0)^{(i)} - \int_0^t \int_U 2\varepsilon \nabla z_N(s)\nabla e_i - \phi'(z_N(s))\phi(z_N(s))|\nabla u_N(s)|^2 e_i \mathrm{d}x\mathrm{d}s \tag{2.11}$$

$$+ \int_0^t \int_U \frac{1 - z_N(s)}{2\varepsilon}e_i \mathrm{d}x\mathrm{d}s$$

for every $1 \le i \le N$. However, due to the strong non linearities at play, *a priori* for every $N$ we only have a local solution in $[0, t_N)$. In order to gain existence in the whole interval $[0, T]$ of Galerkin approximates we will use the following *a priori* estimate, which holds in a slightly more general situation with less regular initial data.

**Proposition 2.5. (Existence of weak approximate solutions in [0,T])**
   Given $(u_0, z_0) \in [L^2(U)]^2$ and a solution $(u_N, z_N)$ of (2.7) in $V_N$, which is a priori well-defined only up to a time $t_N \le T$. It holds

$$\sup_{0 \le t \le t_N} [||u_N(t)||^2_{L^2(U)}] + \int_0^{t_N} ||\nabla u_N(s)||^2_{L^2(U)}\mathrm{d}s \le C$$

and

$$\sup_{0 \le t \le t_N} [||z_N(t)||^2_{L^2(U)}] + \int_0^{t_N} ||\nabla z_N(s)||^2_{L^2(U)}\mathrm{d}s \le C$$

with $C$ a positive constant independent of $N$. By such estimates, we can rule out the possibility of a blow up in finite time for every $N$, thus the same result applies with $T$ in place of $t_N$.

*Proof.* Testing the equation for $u_N$ in (2.7) with $u_N$ itself, we get

$$\frac{1}{2}\frac{\mathrm{d}}{\mathrm{d}t}||u_N(t)||^2_{L^2(U)} = -\int_U (\eta_\varepsilon + \phi(z_N(t))^2)|\nabla u_N(t)|^2\mathrm{d}x - \int_U u_N(t)(u_N(t) - g_N)\mathrm{d}x. \tag{2.12}$$

So

$$\frac{1}{2}\frac{\mathrm{d}}{\mathrm{d}t}||u_N(t)||^2_{L^2(U)} \le ||g_N||_{L^2(U)}||u_N(t)||_{L^2(U)} \le ||g||_{L^2(U)}(1 + ||u_N(t)||^2_{L^2(U)}) \tag{2.13}$$

and one gets uniform boundedness of $||u_N(t)||_{L^2(U)}$ by Gronwall's lemma. Going back to (2.12) one can easily obtain the bound on $||\nabla u_N||_{L^2(0,T;L^2(U))}$ by integrating in time. By orthogonality of the $e_i$ 's we have

$$||u_N(t)||^2_{L^2(U)} = \sum_{i=1}^N [u^{(i)}(t)]^2 \text{ and } ||z_N(t)||^2_{L^2(U)} = \sum_{i=1}^N [z^{(i)}(t)]^2 \qquad (2.14)$$

so we cannot have a blow-up in finite time of the $\{u^{(i)}\}_{i=1}^N$ and we have thus proved existence up to time $T$ for every $T > 0$. The same argument holds if we test the equation for $z_N$ with $z_N$ itself, obtaining

$$\sup_{0 \le t \le T} [||z^N(t)||^2_{L^2(U)}] + 2\varepsilon \int_0^T ||\nabla z_N(s)||^2_{L^2(U)} \mathrm{d}s \le C. \qquad (2.15)$$

$$\square$$

Now that we have existence of solutions of (2.7) in $[0, T]$ for every $N$, let us prove stronger inequalities exploiting the variational characterization of the problem.

**Proposition 2.6.** *(A priori energy estimates)*
Assume $(u_0, z_0) \in [H^1(U)]^2$, then

$$\sup_{t \in [0,T]} [AT_\varepsilon(u_N(t), z_N(t))] + \int_0^T ||\partial_t u_N(t)||^2_{L^2(U)} + ||\partial_t z_N(t)||^2_{L^2(U)} \mathrm{d}t$$

$$= \sup_{t \in [0,T]} [AT_\varepsilon(u_N(t), z_N(t))] + \int_0^T ||\pi_N[\nabla AT_\varepsilon(u_N(t), z_N(t))]||^2_{[L^2(U)]^2} \mathrm{d}t$$

$$\le C,$$

*where $C$ does not depend on $N$. In particular, we have uniform (in $N$) bounds for*

$$(\partial_t u_N, \partial_t z_N) \in [L^2(0,T; L^2(U))]^2 \text{ and } (u_N, z_N) \in [L^\infty(0,T; H^1(U))]^2.$$

*Proof.* Differentiating $AT_\varepsilon(u_N(t), z_N(t))$ in time and exploiting the gradient flow structure of (2.7) we get:

$$\frac{\mathrm{d}}{\mathrm{d}t} AT_\varepsilon(u_N(t), z_N(t)) = -||\pi_N \nabla AT_\varepsilon(u_N(t), z_N(t))||^2_{[L^2(U)]^2} \qquad (2.16)$$

$$= -||\partial_t u_N(t)||^2_{L^2(U)} - ||\partial_t z_N(t)||^2_{L^2(U)}.$$

Integrating the above equality from 0 to $t$:

$$AT_\varepsilon(u_N(t), z_N(t)) + \int_0^t ||\partial_t u_N(s)||^2_{L^2(U)} + ||\partial_t z_N(s)||^2_{L^2(U)} \mathrm{d}s = AT_\varepsilon(\pi_N u_0, \pi_N z_0) \qquad (2.17)$$

$$\le C,$$

with $C$ not depending on $t$. Thus

$$\sup_{t \in [0,T]} [AT_\varepsilon(u_N(t), z_N(t))] + \int_0^T ||\partial_t u_N(t)||^2_{L^2(U)} + ||\partial_t z_N(t)||^2_{L^2(U)} \mathrm{d}t \le C.$$

$$\square$$

**Remark 2.7.** Notice that *a priori* we have no control in $N$ for the quantity

$$\int_U (\eta_\varepsilon + (\pi_N[z_0])^2)|\nabla \pi_N[u_0]|^2 \mathrm{d}x. \tag{2.18}$$

But since we truncated with $|\phi| \leq 2$ we have

$$\int_U (\eta_\varepsilon + \phi(z_0)^2)|\nabla \pi_N[u_0]|^2 \mathrm{d}x \leq \int_U (\eta_\varepsilon + 4)|\nabla u_0|^2 \mathrm{d}x. \tag{2.19}$$

Notice also that the assumption $(u_0, z_0) \in [H^1(U)]^2$ is crucial, otherwise we could have $AT_\varepsilon(u_0, z_0) = +\infty$.

At this point we still cannot prove the weak convergence of the nonlinear parts of the equation, namely $\mathrm{div}((\eta_\varepsilon + \phi(z_N)^2)\nabla u_N)$ and $\phi'(z_N)\phi(z_N)|\nabla u|^2$. Stronger estimates are needed. The next Proposition will be the crucial one, providing uniform $[L^2(0,T;H^2(U))]^2$ boundedness of $(u_N, z_N)$. Before that we state a version of the Gagliardo–Nirenberg interpolation inequality which will be used in the proof.

**Lemma 2.8** (Gagliardo–Nirenberg, [11]). *Let $\Omega$ be a bounded $C^2$ domain in $\mathbb{R}^n$. Assume that $1 \leq q, r \leq +\infty$ such that*

$$\frac{1}{p} = \frac{j}{n} + \theta\left(\frac{1}{r} - \frac{m}{n}\right) + \frac{1-\theta}{q}, \ \text{with } \frac{j}{m} \leq \theta \leq 1. \tag{2.20}$$

*Then*

$$||D^j u||_{L^p(\Omega)} \leq C||D^m u||^\theta_{L^r(\Omega)}||u||^{1-\theta}_{L^q(\Omega)} + C||u||_{L^s(\Omega)},$$

*where $s > 0$ is arbitrary and the constant $C$ depends on $\Omega, q, r, j, m, n, \theta, s$.*

**Remark 2.9.** Sharper versions of Lemma 2.8 are available in the literature, we just stated a version which is sufficient for our purposes.

**Proposition 2.10.** *(Uniform $L^2(0,T;H^2(U))$ estimates)*
Let $U \subset \mathbb{R}^2$ be a bounded $C^2$ domain and let $(u_0, z_0) \in [H^1(U)]^2$. Consider solutions $(u_N, z_N)$ of (2.7). Then it holds that

$$\int_0^T ||u_N(t)||^2_{H^2(U)} + ||z_N(t)||^2_{H^2(U)} + ||\nabla u_N(t)||^4_{L^4(U)} + ||\nabla z_N(t)||^4_{L^4(U)} \mathrm{d}t \leq C \tag{2.21}$$

*for some $C > 0$ independent of $N$.*

*Proof.* First of all, we notice that having already proved in Proposition 2.6 a bound for the $L^\infty(0,T;H^1(U))$ norm of both $u_N$ and $z_N$, to prove (2.21) it will be sufficient to show that

$$\int_0^T ||\Delta u_N(t)||^2_{L^2(U)} + ||\Delta z_N(t)||^2_{L^2(U)} + ||\nabla u_N(t)||^4_{L^4(U)} + ||\nabla z_N(t)||^4_{L^4(U)} \mathrm{d}t \leq C.$$

We divide the proof in steps to improve readability.

- **STEP 1: Control of the $L^2(0, T; H^2(U))$ norm with the $L^4(0, T; L^4(U))$ norm of the gradients**
  We want to prove an estimate like

$$\sup_{0 \leq t \leq T} [AT_\varepsilon(u_N(t), z_N(t))] + \int_0^T ||\Delta u_N(t)||^2_{L^2(U)} + ||\Delta z_N(t)||^2_{L^2(U)} \mathrm{d}t \tag{2.22}$$
$$\lesssim 1 + \int_0^T ||\nabla u_N(t)||^4_{L^4(U)} + ||\nabla z_N(t)||^4_{L^4(U)} \mathrm{d}t.$$

The idea is to expand the energy equality (2.17) obtained in Proposition 2.6 with $||\pi_N \nabla AT_\varepsilon(u_N, z_N)||^2_{L^2(U)}$. For the sake of readability we omit writing time and space dependence, abbreviate $\phi(z_N)$ as $\phi$ and omit the subscripts too:

$$\int_0^T ||\nabla AT_\varepsilon(u, z)||^2_{L^2(U)^2} \mathrm{d}t = \int_0^T ||\partial_u AT_\varepsilon(u, z)||^2_{L^2(U)} + ||\partial_z AT_\varepsilon(u, z)||^2_{L^2(U)} \mathrm{d}t \tag{2.23}$$
$$= \int_0^T \int_U [(\eta + \phi^2)\Delta u + 2\phi\phi'\nabla z \nabla u - (u - g)]^2 + \left[2\varepsilon\Delta z - \phi\phi'|\nabla u|^2 + \frac{1-z}{2\varepsilon}\right]^2 \mathrm{d}x \mathrm{d}t$$
$$\geq \int_0^T \int_U \frac{(\eta + \phi^2)^2(\Delta u)^2}{3} - 12\phi^2(\phi')^2(\nabla z \nabla u)^2 - 3(u - g)^2 \mathrm{d}x \mathrm{d}t$$
$$+ \int_0^T \int_U \frac{4\varepsilon^2(\Delta z)^2}{3} - 3\phi^2(\phi')^2|\nabla u|^4 - \frac{3(1-z)^2}{4\varepsilon^2} \mathrm{d}x \mathrm{d}t$$
$$\gtrsim \int_0^T ||\Delta u||^2_{L^2(U)} + ||\Delta z||^2_{L^2(U)} \mathrm{d}t - \int_0^T ||\nabla u||^4_{L^4(U)} + ||\nabla z||^4_{L^4(U)} \mathrm{d}t$$
$$- ||u||_{L^2(0,T;L^2(U))} - ||z||_{L^2(0,T;L^2(U))} - ||g||_{L^2(U)} - 1,$$

where we have used the inequalities $(a + b + c)^2 \geq a^2/3 - 3b^2 - 3c^2$ and $ab \leq a^2/2 + b^2/2$. The former one in particular follows by noticing that

$$(a + b + c)^2 - \frac{a^2}{3} + 3b^2 + 3c^2 = \left(\frac{1}{\sqrt{3}}a + \sqrt{3}b\right)^2 + \left(\frac{1}{\sqrt{3}}a + \sqrt{3}c\right)^2 + (b + c)^2. \tag{2.24}$$

Using Proposition 2.6, the uniform bounds on $||u_N||_{L^2(0,T;L^2(U))}, ||z_N||_{L^2(0,T;L^2(U))}$ and (2.23) we can easily deduce (2.22).

- **STEP 2: Control of (2.21) with the $L^4(0, T; L^4(U))$ norm of $\nabla u_N$**
  We now reduce proving (2.21) to estimating the $L^4(0, T; L^4(U))$ norm of $\nabla u_N$. Testing the equation in $z_N$ with $-\Delta z_N$ in (2.7) we get

$$\frac{1}{2}\frac{\mathrm{d}t}{\mathrm{d}t}||\nabla z_N(t)||^2_{L^2(U)} = -2\varepsilon||\Delta z_N(t)||^2_{L^2(U)} \tag{2.25}$$
$$+ \int_U \phi'(z_N(t))\phi(z_N(t))|\nabla u_N(t)|^2 \Delta z_N(t) \mathrm{d}x - \int_U \frac{1 - z_N(t)}{2\varepsilon}\Delta z_N(t) \mathrm{d}x.$$

Using that $ab \leq (\delta/2)a^2 + (1/2\delta)b^2$ for a suitable $\delta$ (for example, $\delta = \varepsilon$ works) we have:

$$\frac{1}{2}\frac{\mathrm{d}}{\mathrm{d}t}||\nabla z_N(t)||^2_{L^2(U)} \leq -2\varepsilon||\Delta z_N(t)||^2_{L^2(U)} + \frac{\varepsilon}{2}||\Delta z_N(t)||^2_{L^2(U)} \tag{2.26}$$
$$+ \frac{1}{2\varepsilon}\int_U (\phi'(z_N(t))\phi(z_N(t)))^2|\nabla u_N(t)|^4 \mathrm{d}x + \frac{\varepsilon}{2}||\Delta z_N(t)||^2_{L^2(U)}$$
$$+ \frac{1}{2\varepsilon}\int_U \left(\frac{1 - z_N(t)}{2\varepsilon}\right)^2 \mathrm{d}x,$$

and integrating this inequality in time we have

$$\int_0^T ||z_N(t)||_{H^2(U)}^2 \mathrm{d}t \lesssim 1 + \int_0^T ||\nabla u_N(t)||_{L^4(U)}^4 \mathrm{d}t, \tag{2.27}$$

where we remark that the multiplicative constants in the inequality depend heavily on $\varepsilon$ (in particular, they go to $+\infty$ as $\varepsilon \to 0$). Moreover by Lemma 2.8 with $j = 0, m = 1, p = 4, r = q = s = 2, \theta = 1/2$:

$$||\nabla z_N(t)||_{L^4(U)} \lesssim ||\nabla z_N(t)||_{L^2(U)} + ||\nabla z_N(t)||_{L^2(U)}^{1/2} ||\nabla^2 z_N(t)||_{L^2(U)}^{1/2} \tag{2.28}$$
$$\lesssim 1 + ||\nabla^2 z_N(t)||_{L^2(U)}^{1/2}$$

thanks to the $L^\infty(0, T; H^1(U))$ estimates on $z_N$. Using that for $a, b > 0$ it holds $(a + b)^4 \leq 8(a^4 + b^4)$ then:

$$||\nabla z_N(t)||_{L^4(U)}^4 \lesssim 1 + ||\nabla^2 z_N(t)||_{L^2(U)}^2. \tag{2.29}$$

Combining (2.22), (2.27) and (2.29) we finally get

$$\int_0^T ||u_N(t)||_{H^2(U)}^2 + ||z_N(t)||_{H^2(U)}^2 + ||\nabla u_N(t)||_{L^4(U)}^4 + ||\nabla z_N(t)||_{L^4(U)}^4 \mathrm{d}t \tag{2.30}$$
$$\lesssim 1 + \int_0^T ||\nabla u_N(t)||_{L^4(U)}^4 + ||\nabla z_N(t)||_{L^4(U)}^4 \mathrm{d}t \lesssim 1 + \int_0^T ||\nabla u_N(t)||_{L^4(U)}^4 \mathrm{d}t.$$

- **STEP 3: Control of $||\nabla u_N||_{L^4(0,T;L^4(U))}$ and conclusion**
  Having proved (2.30) we can conclude the proof if we can find a uniform estimate for $||u_N||_{L^4(0,T;L^4(U))}$. Considering the time $t$ fixed we focus on the first equation of (2.7) and we consider $u_N(t)$ the solution of:

$$\begin{cases} -\mathrm{div}((\eta_\varepsilon + \phi(z_N(t))^2)\nabla u_N(t)) = f_N(t) & \text{in } U \\ \partial_n u_N(t) = 0 & \text{in } \partial U \end{cases} \tag{2.31}$$

where $f_N(t) = -\partial_t u_N(t) - u_N(t) + g_N$. By classical elliptic regularity, we have

$$||\nabla^2 u_N(t)||_{L^2(U)} \leq C||f_N(t)||_{L^2(U)},$$

so

$$\int_0^T ||\nabla^2 u_N(t)||_{L^2(U)}^2 \mathrm{d}t \lesssim \int_0^T ||\partial_t u_N(t)||_{L^2(U)}^2 + ||u_N(t)||_{L^2(U)}^2 + ||g_N||_{L^2(U)}^2 \mathrm{d}t \leq C \tag{2.32}$$

by Proposition 2.6. Applying Lemma 2.8 with $j = 0, m = 1, p = 4, r = 2, q = 3, \theta = 1/4$ and $s = 2$ to $\nabla u_N(t)$ we have:

$$||\nabla u_N(t)||_{L^4(U)} \lesssim 1 + ||\nabla^2 u_N(t)||_{L^2(U)}^{1/4} ||\nabla u_N(t)||_{L^3(U)}^{3/4}, \tag{2.33}$$

so, using (2.33) and (2.32):

$$\int_0^T ||\nabla u_N(t)||_{L^4(U)}^4 \mathrm{d}t \lesssim 1 + \int_0^T ||\nabla^2 u_N(t)||_{L^2(U)} ||\nabla u_N(t)||_{L^3(U)}^3 \mathrm{d}t \tag{2.34}$$

$$\lesssim 1 + \left( \int_0^T ||\nabla^2 u_N(t)||_{L^2(U)}^2 \mathrm{d}t \right)^{1/2} \left( \int_0^T ||\nabla u_N(t)||_{L^3(U)}^6 \mathrm{d}t \right)^{1/2}$$

$$\lesssim 1 + \left( \int_0^T ||\nabla u_N(t)||_{L^3(U)}^6 \mathrm{d}t \right)^{1/2}.$$

To bound this last term, we use once again Lemma 2.8 with $j = 0, m = 1, p = 3, r = q = 2, \theta = 1/3$ and $s = 2$ to get:

$$||\nabla u_N(t)||_{L^3(U)} \lesssim 1 + ||\nabla^2 u_N(t)||_{L^2(U)}^{1/3} ||\nabla u_N(t)||_{L^2(U)}^{2/3}.$$

So:

$$\int_0^T ||\nabla u_N(t)||_{L^3(U)}^6 \mathrm{d}t \lesssim 1 + \int_0^T ||\nabla^2 u_N(t)||_{L^2(U)}^2 ||\nabla u_N(t)||_{L^2(U)}^4 \mathrm{d}t \leq C,$$

because $u_N \in L^\infty(0,T;H^1(U))$ and because of (2.32).

$\square$

**Remark 2.11.** It seems to us that being in dimension 2 is crucial for Proposition 2.10 to hold. Assume, for example, to be in dimension 3. Then if we wanted to obtain the equivalent of (2.29) we would have to use Lemma 2.8 as:

$$||\nabla z_N(t)||_{L^4(U)} \lesssim 1 + ||\nabla z_N(t)||_{L^q(U)}^{1-\theta} ||\nabla^2 z_N(t)||_{L^2(U)}^\theta$$

with $q \in [1, +\infty]$ and $\theta \in [0,1]$ such that (2.20) holds, namely:

$$\frac{1}{4} = \frac{\theta}{6} + \frac{1-\theta}{q}. \tag{2.35}$$

Since we want to bound $||\nabla z_N(t)||_{L^q(U)}$, by Sobolev embedding we can hope to do so only if $q \leq 6 = 2^*$, and it is not difficult to see that (2.35) can hold only if $q \in [1,4]$. Let us consider two cases:

- $1 \leq q \leq 2$
  By (2.35) we have $\theta = (12 - 3q)/(12 - 2q)$ and we can use the boundedness of $U$ and the $L^\infty$ estimate on $||\nabla z_N(t)||_{L^2(U)}$ to get, by Hölder's inequality:

  $$||\nabla z_N(t)||_{L^4(U)} \lesssim 1 + ||\nabla z_N(t)||_{L^4(U)}^\theta.$$

The value of $\theta$ as a function of $q$ is decreasing in $[1,4]$, so at best when $q = 2$ we get

$$||\nabla z_N(t)||_{L^4(U)} \lesssim 1 + ||\nabla z_N(t)||_{L^4(U)}^{3/4} \implies ||\nabla z_N(t)||_{L^4(U)}^4 \lesssim 1 + ||\nabla z_N(t)||_{L^4(U)}^3 \tag{2.36}$$

which is not enough to conclude.

- $2 < q \le 4$

  In this case it seems to us that we would be forced to use the Sobolev inequality to bound $||\nabla z_N(t)||_{L^q(U)}$, thus for every $\theta$ we would get

  $$||\nabla z_N(t)||_{L^4(U)} \lesssim 1 + ||\nabla z_N(t)||_{L^q(U)}^{1-\theta} ||\nabla^2 z_N(t)||_{L^2(U)}^\theta \lesssim 1 + ||\nabla^2 z_N(t)||_{L^2(U)}, \tag{2.37}$$

  which is clearly not enough to conclude.

Of course this is not a counterexample, but it seems to suggest that if we hope to prove regularity in higher dimensions, different methods might be needed.

Before proving the main result, let us prove the following Lemma:

**Lemma 2.12.** *Let $X$ be a subset of $\mathbb{R}^n$ of finite measure and let $\{f_j\}$ be a sequence of functions in $L^p(X)$ with $1 < p < +\infty$ for which $\exists M > 0$ such that $||f_j||_{L^p(X)} \le M$ for every $j$ and $f_j(x) \to f(x)$ $\mathcal{L}^n$-almost everywhere in $X$. Then*

$$f_j \rightharpoonup f \ in \ L^p(X).$$

*Proof.* The proof of a stronger version of this lemma can be found in [12], but we report it here since it is very short. By Fatou's lemma, it holds

$$||f||_{L^p(X)} \le \liminf_{j \to +\infty} ||f_j||_{L^p(X)} \le M.$$

Fix $g \in L^q(X)$ with $1/p + 1/q = 1$ and fix $\varepsilon > 0$. By Egoroff's theorem there exists a subset $E \subset X$ such that $f_j$ converges to $f$ uniformly in $E$ and $||g||_{L^q(X \setminus E)} \le \varepsilon/(2M)$. Then:

$$\left| \int_X f_j g \mathrm{d}x - \int_X f g \mathrm{d}x \right| \le \int_X |f_j - f| |g| \mathrm{d}x = \int_{X \setminus E} |f - f_j| |g| \mathrm{d}x + \int_E |f - f_j| |g| \mathrm{d}x$$
$$\le ||f_j - f||_{L^p(X \setminus E)} ||g||_{L^q(X \setminus E)} + ||f_j - f||_{L^p(E)} ||g||_{L^q(E)}.$$

Now:

$$||f - f_j||_{L^p(X \setminus E)}^p \le 2^{p-1} \left( \int_{X \setminus E} |f|^p \mathrm{d}x + \int_{X \setminus E} |f_j|^p \mathrm{d}x \right) \le 2^p M^p,$$

so

$$||f - f_j||_{L^p(X \setminus E)} \le 2M.$$

Moreover, since $X$ has finite measure:

$$||f_j - f||_{L^p(E)} \le \mathcal{L}^n(X)^{1/p} \sup_{x \in E} |f_j(x) - f(x)| \to 0 \text{ as } j \to +\infty.$$

Thus we have that:

$$\limsup_j \left| \int_X f_j g \mathrm{d}x - \int_X f g \mathrm{d}x \right| \le 2M ||g||_{L^q(X \setminus E)} \le \varepsilon,$$

and we conclude by the arbitrarity of $\varepsilon$. $\qquad\square$

We are now ready to prove Theorem 2.2

*Proof of Theorem 2.2.* The uniqueness of the solution follows exactly as in Step 3 of [8], Theorem 2.3 and we will give it for granted, not repeating the proof here. Up to a subsequence we will not rename, given the uniform bounds we proved on Proposition 2.6 and Proposition 2.10 we consider $(u, z)$ such that:

$$
\begin{cases}
(u_N, z_N) \rightharpoonup (u, z) & \text{in } [L^2(0, T; H^2(U))]^2 \\
(\partial_t u_N, \partial_t z_N) \rightharpoonup (\partial_t u, \partial_t z) & \text{in } [L^2(0, T; L^2(U))]^2 \\
(u_N, z_N) \to (u, z) & \text{in } [C(0, T; L^2(U))]^2 \\
(u_N, z_N) \to (u, z) & \text{in } [L^2(0, T; H^1(U))]^2 \\
(u_N, z_N, \nabla u_N, \nabla z_N) \to (u, z, \nabla u, \nabla z) & \text{pointwise almost everywhere,}
\end{cases}
\tag{2.38}
$$

where the compact embeddings in $C(0, T; L^2(U))$ and $L^2(0, T; H^1(U))$ are obtained by applying the Aubin–Lions lemma (see [13], [14], [15]). Notice that we also know that

$$
(\partial_t u_N, \partial_t z_N) \rightharpoonup (\partial_t u, \partial_t z) \text{ in } [L^2(0, T, L^2(U))]^2
\tag{2.39}
$$

(see [16], Chap. 7, Probl. 5) and

$$
|\nabla u_N|^2 \rightharpoonup |\nabla u|^2 \text{ in } L^2(0, T; L^2(U)).
\tag{2.40}
$$

The weak convergence of $|\nabla u_N|^2$ to $|\nabla u|$ follows by identifying $L^2(0, T; L^2(U))$ with $L^2([0, T] \times U)$ and applying Lemma 2.12. By (2.38) we know

$$
(u, z) \in [L^2(0, T; H^2(U)) \cap L^\infty(0, T; H^1(U)) \cap H^1(0, T; L^2(U))]^2,
$$

so to conclude we only need to show that $(u, z)$ is a solution of (1.4), and thanks to Proposition 2.4 it will be enough to prove that $(u, z)$ solves (2.1). Let $\psi \in V_M = \text{Span}\{e_1, \ldots, e_M\}$ be a test function for (2.7) with $N > M$, so it holds:

$$
\underbrace{\int_U u_N(t)\psi(x)\mathrm{d}x}_{\textcircled{1}} = \int_U \pi_N[u_0]\psi(x)\mathrm{d}x - \underbrace{\int_0^t \int_U (\eta_\varepsilon + \phi(z_N)^2)\nabla u_N \nabla \psi \mathrm{d}x\mathrm{d}s}_{\textcircled{2}}
$$
$$
- \int_0^t \int_U (u_N - g_N)\psi\mathrm{d}x\mathrm{d}s
\tag{2.41}
$$

and

$$
\underbrace{\int_U z_N(t)\psi(x)\mathrm{d}x}_{\textcircled{3}} = \int_U \pi_N[z_0]\psi(x)\mathrm{d}x - 2\varepsilon \int_0^t \int_U \nabla z_N \nabla \psi \mathrm{d}x\mathrm{d}s
$$
$$
- \underbrace{\int_0^t \int_U \phi'(z_N)\phi(z_N)|\nabla u_N|^2 \psi \mathrm{d}x\mathrm{d}s}_{\textcircled{4}} + \int_0^t \int_U \frac{1 - z_N}{2\varepsilon}\psi\mathrm{d}x\mathrm{d}s.
\tag{2.42}
$$

We want to show we can pass to the limit in $N$ in every highlighted term, since for the others it is trivial by weak convergence.

- As for ① and ③, we can pass to the limit thanks to the compactness in $C(0, T; L^2(U))$.
- For ②, we have (by dominated convergence) strong convergence in
  $L^2(0, T; L^2(U))$ of $(\eta_\varepsilon + \phi(z_N)^2)$, and weak convergence of $\nabla u_N$. So their product weakly converges since $(\eta_\varepsilon + \phi(z_N)^2)$ is uniformly bounded and we can pass to the limit.
- We have that $\phi'(z_N)\phi(z_N) \to \phi'(z)\phi(z)$ in $L^2(0, T; L^2(U))$ and it is uniformly bounded, so the product with $|\nabla u_N|^2$ weakly converges and we can pass to the limit in ④.

It only remains to prove that $(u, z)$ satisfies the homogeneous Neumann boundary conditions of (1.4). Let $\psi \in L^2(0, T; H^1(U))$, then

$$0 = \int_0^T \int_{\partial U} \psi \partial_n u_N \mathrm{d}\mathcal{H}^1(x) \mathrm{d}t = \int_0^T \int_U \nabla u_N \nabla \psi \mathrm{d}x \mathrm{d}t + \int_0^T \int_U \Delta u_N \psi \mathrm{d}x \mathrm{d}t.$$

Since $u_N \rightharpoonup u$ in $L^2(0, T; H^2(U))$, we can pass to the limit in both integrals and we get:

$$0 = \int_0^T \int_U \nabla u \nabla \psi \mathrm{d}x \mathrm{d}t + \int_0^T \int_U \Delta u \psi \mathrm{d}x \mathrm{d}t = \int_0^T \int_{\partial U} \psi \partial_n u \mathrm{d}\mathcal{H}^1(x) \mathrm{d}.$$

By the arbitrarity of $\psi$, we have that $\partial_n u = 0$ almost everywhere on $(0, T) \times \partial U$, and the same procedure works for $z_N$. The proof is complete. $\qquad\square$

## Acknowledgements

The author wishes to thank an anonymous reviewer whose help and suggestions have helped to simplify some proofs and improve the readability of the paper. The author is a member of the INdAM group GNAMPA.

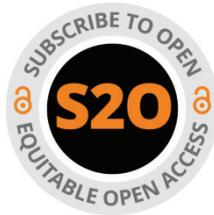

**Please help to maintain this journal in open access!**

This journal is currently published in open access under the Subscribe to Open model (S2O). We are thankful to our subscribers and supporters for making it possible to publish this journal in open access in the current year, free of charge for authors and readers.

Check with your library that it subscribes to the journal, or consider making a personal donation to the S2O programme by contacting `subscribers@edpsciences.org`.

More information, including a list of supporters and financial transparency reports, is available at https://edpsciences.org/en/subscribe-to-open-s2o.